\documentclass[12pt]{article}
\usepackage{amsbsy,amsfonts,amsmath,amssymb,amsthm,enumerate,graphicx,bbm}

\numberwithin{equation}{section}
\theoremstyle{plain}

\newtheorem{theorem}{Theorem}[section]
\newtheorem{corollary}[theorem]{Corollary}

\newcommand{\N}{\mathbb{N}}

\newcommand{\dc}{d_{\rm c}}

\newcommand{\lbeq}[1]{\label{eq:#1}}
\newcommand{\laplace}{\bar\Delta}

\newcommand{\mc}{m_{\rm c}}

\newcommand{\mR}{{\mathbb R}}
\newcommand{\mZ}{{\mathbb Z}}
\newcommand{\nablone}{\nabla_{\!1}}
\newcommand{\nn}{\nonumber}

\newcommand{\pc}{p_{\rm c}}

\newcommand{\Rd}{\mR^d}

\newcommand{\refeq}[1]{(\ref{eq:#1})}
\newcommand{\sss}{\scriptscriptstyle}

\newcommand{\Zd}{\mZ^d}
\newcommand{\Zp}{\mZ_+}

\newcommand{\op}{{\sss\rm OP}}
\newcommand{\rw}{{\sss\rm RW}}
\newcommand{\saw}{{\sss\rm SAW}}

\usepackage[usenames]{color}

\title{Large-time asymptotics of the gyration radius for long-range
statistical-mechanical models\footnote{Submitted to
\emph{RIMS Kokyuroku Bessatsu}.}}

\author{
Akira~Sakai\footnote{Creative Research Institution SOUSEI,
Hokkaido University, Japan.  {\tt sakai@cris.hokudai.ac.jp}}
}

\begin{document}
\maketitle

\begin{abstract}
The aim of this short article is to convey the basic idea of the original
paper \cite{csIII}, without going into too much detail, about how to derive
sharp asymptotics of the gyration radius for random walk, self-avoiding walk
and oriented percolation above the model-dependent upper critical dimension.
\end{abstract}

\section{Introduction}
Let $D$ be the $\Zd$-symmetric 1-step distribution for random walk (RW) and
define the RW 2-point function as
\begin{align}\lbeq{rw2pt}
\varphi_t^\rw(x)=\sum_{\substack{\omega:o\to x\\ |\omega|=t}}\prod_{s=1}^t
 D(\omega_s-\omega_{s-1})\qquad(x\in\Zd,~t\in\Zp).
\end{align}
We also consider self-avoiding walk (SAW) and oriented percolation (OP) that
are both generated by $D$.  The SAW 2-point function is defined as
\begin{align}
\varphi_t^\saw(x)=\sum_{\substack{\omega:o\to x\\ |\omega|=t}}\prod_{s=1}^t
 D(\omega_s-\omega_{s-1})\prod_{0\le i<j\le t}(1-\delta_{\omega_i,\omega_j}),
\end{align}
where the indicator $\prod_{0\le i<j\le t}(1-\delta_{\omega_i,\omega_j})$,
which is absent in \refeq{rw2pt}, is 1 if and only if $\omega$ does not
intersect to itself, hence accounting for the self-avoidance constraint.
Oriented percolation is a model for random media in space-time $\Zd\times\Zp$.
A bond is an ordered pair of vertices in $\Zd\times\Zp$, and each bond
$((u,t),(v,t+1))$ is either occupied or vacant with probability $pD(v-u)$ and
$1-pD(v-u)$, respectively, independently of the other bonds.  The parameter
$p$ equals the expected number of occupied bonds per vertex, and it is known
that there is a phase transition at $p=\pc$.  We say that $(x,s)$ is connected
to $(y,t)$ if either $(x,s)=(y,t)$ or there is a time-increasing sequence of
occupied bonds from $(x,s)$ to $(y,t)$.  The OP 2-point function
$\varphi_t^\op(x)$ is then defined as the probability that the origin $(o,0)$
is connected to $(x,t)$.

The models are said to be finite-range if $D$ is supported on a finite set of
$\Zd$.  The main property of a finite-range $D$ is the existence of the
variance $\sigma^2\equiv\sum_{x\in\Zd}|x|^2D(x)$, and because of this,
investigation of finite-range models is relatively easier.  The situation is
basically the same for $D$ that decays faster than any polynomials, such as an
exponentially decaying $D$.  However, if $D(x)\approx|x|^{-d-\alpha}$ for
large $|x|$, then the existence of the variance depends on $\alpha>0$ and
therefore we cannot always expect that the same results for finite-range
models also hold for this long-range models with index $\alpha$.  For
example, take the gyration radius of order $r\in(0,\alpha)$, which is defined
as
\begin{align}\lbeq{gyradius}
\xi_t^{\sss(r)}=\bigg(\frac{\sum_{x\in\Zd}|x|^r\varphi_t(x)}{\sum_{x\in\Zd}
 \varphi_t(x)}\bigg)^{1/r}.
\end{align}
The gyration radius represents a typical end-to-end distance of a linear
structure of length $t$ or a typical spatial size of a cluster at time $t$.
It may be natural to guess, at least for random walk, that
$\xi_t^{\sss(r)}=O(\sqrt t)$ if $\alpha>2$ and
$\xi_t^{\sss(r)}=O(t^{1/\alpha})$ if $\alpha<2$, for every real
$r\in(0,\alpha)$.  As we state shortly, we have proved affirmative results
\cite{csIII} for random walk in any dimension and for self-avoiding walk and
critical/subcritical oriented percolation above the common upper-critical
dimension $\dc\equiv2(\alpha\wedge2)$.

More precisely, we assume the following properties of $D$.  Given an
$L\in[1,\infty)$, we suppose that $D(x)\propto|x/L|^{-d-\alpha}$
for large $|x|$ such that its Fourier transform
$\hat D(k)\equiv\sum_{x\in\Zd}e^{ik\cdot x}D(x)$ exhibits the $k\to0$
asymptotics
\begin{align}\lbeq{D-asympt}
1-\hat D(k)=v_\alpha|k|^{\alpha\wedge2}\times
 \begin{cases}
 1+O((L|k|)^\epsilon)&(\alpha\ne2),\\
 \log\frac1{L|k|}+O(1)\quad&(\alpha=2),
 \end{cases}
\end{align}
for some $v_\alpha=O(L^{\alpha\wedge2})$ and $\epsilon>0$.  If $\alpha>2$ (or
$D$ is finite-range), then $v_\alpha=\frac1{2d}\sigma^2$.  An example that
satisfies the above properties is the long-range Kac potential
\begin{align}\lbeq{kac}
D(x)=\frac{h(y/L)}{\sum_{y\in\Zd}h(y/L)}\qquad(x\in\Zd),
\end{align}
defined by the rotation-invariant function
\begin{align}\lbeq{h-def}
h(x)=\frac{1+O\big((|x|\vee1)^{-\rho}\big)}{(|x|\vee1)^{d+\alpha}}\qquad
 (x\in\Rd),
\end{align}
for some $\rho>\epsilon$ (cf., \cite{csIII}).  Under this assumption, we have
proved the following sharp asymptotics of a variant of the gyration radius:

\begin{theorem}[\cite{csIII}]\label{thm:main}
For random walk in any dimension with any $L$, and for self-avoiding walk and
critical/subcritical oriented percolation for $d>\dc$
with $L\gg1$, there is a model-dependent constant $C_\alpha=1+O(L^{-d})$
($C_\alpha\equiv1$ for random walk) such that, for every $r\in(0,\alpha)$,
\begin{align}\lbeq{main2}
\frac{\sum_{x\in\Zd}|x_1|^r\varphi_t(x)}{\sum_{x\in\Zd}\varphi_t(x)}
&\underset{t\uparrow\infty}\sim\frac{2\sin\frac{r\pi}{\alpha\vee2}}{(\alpha
 \wedge2)\sin\frac{r\pi}\alpha}\,\frac{\Gamma(r+1)}{\Gamma(\frac{r}{\alpha
 \wedge2}+1)}\times
 \begin{cases}
 (C_\alpha v_\alpha t)^{\frac{r}{\alpha\wedge2}}&(\alpha\ne2),\\
 (C_2v_2t\log\sqrt{t})^{r/2}&(\alpha=2),
 \end{cases}
\end{align}
where $x_1$ is the first coordinate of $x\in\Zd$.
\end{theorem}

We should emphasize that, except for the actual value of $C_\alpha$, the
expression \refeq{main2} is universal.  The result also holds for finite-range
models, for which $\alpha$ is considered to be infinity.  As far as we notice,
even for random walk, the sharp asymptotic expression \refeq{main2} for all real
$r\in(0,\alpha)$ is new.

Using $|x_1|^r\le|x|^r\le d^{r/2}\sum_{j=1}^d|x_j|^r$ and the $\Zd$-symmetry of
the models, we can conclude the following:

\begin{corollary}[\cite{csIII}]
Under the same condition as in Theorem~\ref{thm:main},
\begin{align}\lbeq{conjecture}
\xi_t^{\scriptscriptstyle(r)}=
 \begin{cases}
 O(t^{\frac1{\alpha\wedge2}})&(\alpha\ne2),\\
 O(\sqrt{t\log t})\quad&(\alpha=2).
 \end{cases}
\end{align}
for every $r\in(0,\alpha)$.
\end{corollary}

In his recent work \cite{h??}, Heydenreich proved \refeq{conjecture} for
self-avoiding walk, but only for small $r<\alpha\wedge2$, with no attempt to
identify the proportional constant.  Our results are somewhat stronger,
because we have derived the exact expression for the proportional constant in
\refeq{main2} (also clarifying its model-dependence) and proved
\refeq{conjecture} for all $r<\alpha$.

\section{Sketch proof for random walk}\label{s:rw}
In this section, we restrict our attention to random walk, which is obviously
simpler than the other two models, and explain the framework of the proof of
Theorem~\ref{thm:main}.

First we consider the generating function (= the Fourier-Laplace transform)
of the 2-point function.  Recall that $\varphi_t^\rw(x)$ satisfies the
convolution equation
\begin{align}\lbeq{rw-exp-x}
\varphi_t^\rw(x)=\delta_{t,0}\delta_{x,o}+(D*\varphi_{t-1}^\rw)(x)\equiv
 \delta_{t,0}\delta_{x,o}+\sum_{y\in\Zd}D(y)\,\varphi_{t-1}^\rw(x-y),
\end{align}
where we regard $(D*\varphi_{t-1}^\rw)(x)$ for $t\le0$ as zero.  Taking the 
Fourier-Laplace transform of both sides, we obtain that, for
$k\in[-\pi,\pi]^d$ and $m\in[0,\mc^\rw)$,
\begin{align}\lbeq{rw-exp}
\hat\varphi_m^\rw(k)\equiv\sum_{t\in\Zp}m^t\sum_{x\in\Zd}e^{ik\cdot x}
 \varphi_t^\rw(x)=1+m\hat D(k)\,\hat\varphi_m^\rw(k),
\end{align}
where $\mc^\rw\equiv1$ is the radius of convergence for the sequence
$\big\{\sum_{x\in\Zd}\varphi_t^\rw(x)\big\}_{t\in\Zp}$.  To see this in a
different way, take $k=0$ in \refeq{rw-exp} so that
\begin{align}\lbeq{rw-initcond}
\hat\varphi_m^\rw(0)=1+m\hat\varphi_m^\rw(0)=\frac1{1-m}.
\end{align}
The expansion of the right-hand side is $\sum_{t\in\Zp}m^t$ and the
coefficient of $m^t$ is exactly 1 ($\equiv\sum_{x\in\Zd}\varphi_t^\rw(x)$)
for every $t\in\Zp$.

Next we differentiate $\hat\varphi_m^\rw(k)$ with respect to $k_1$ (= the
first coordinate of $k$) to yield the generating function of the sequence
$\big\{\sum_{x\in\Zd}|x_1|^r\varphi_t^\rw(x)\big\}_{t\in\Zp}$.  For example,
if $r=2j$ with $j\in\N$ (hence $\alpha>2$), then
\begin{align}\lbeq{nablone-def}
\nablone^{2j}\hat\varphi_m^\rw(0)\equiv\frac{\partial^{2j}}{\partial k_1^{2j}}
 \hat\varphi_m^\rw(k)\bigg|_{k=0}=(-1)^j\sum_{t\in\Zp}m^t\sum_{x\in\Zd}x_1^{2j}
 \varphi_t^\rw(x).
\end{align}
On the other hand, by differentiating \refeq{rw-exp} and using the
$\Zd$-symmetry of the model,
\begin{align}
\nablone^{2j}\hat\varphi_m^\rw(0)&=m\nablone^{2j}\hat\varphi_m^\rw(0)+m\sum_{l
 =1}^j\binom{2j}{2l}\nablone^{2l}\hat D(0)\,\nablone^{2(j-l)}\hat\varphi_m^\rw
 (0)\nn\\
&=\frac{m}{1-m}\sum_{l=1}^j\binom{2j}{2l}\nablone^{2l}\hat D(0)\,\nablone^{2(j
 -l)}\hat\varphi_m^\rw(0).
\end{align}
Solving this recursion by induction under the initial condition
\refeq{rw-initcond}, we obtain (see \cite{csIII} for more details)
\begin{align}
\nablone^{2j}\hat\varphi_m^\rw(0)&=\binom{2j}2\frac{m\nablone^2\hat D(0)}{1-m}
 \,\nablone^{2(j-1)}\hat\varphi_m^\rw(0)+O\big((1-m)^{-j}\big)\nn\\
&=\binom{2j}2\binom{2(j-1)}2\bigg(\frac{m\nablone^2\hat D(0)}{1-m}\bigg)^2
 \nablone^{2(j-2)}\hat\varphi_m^\rw(0)+O\big((1-m)^{-j}\big)\nn\\
&~\;\vdots\nn\\
&=\prod_{l=1}^j\binom{2l}2\bigg(\frac{m\nablone^2\hat D(0)}{1-m}\bigg)^j\hat
 \varphi_m^\rw(0)+O\big((1-m)^{-j}\big)\nn\\
&=\frac{(2j)!}{2^j}\frac{\big(m\nablone^2\hat D(0)\big)^j}{(1-m)^{j+1}}+O\big(
 (1-m)^{-j}\big).
\end{align}
Comparing this with \refeq{nablone-def} and using
$v_\alpha\equiv\frac1{2d}\sigma^2=\frac{-1}2\nablone^2\hat D(0)$ for
$\alpha>2$, we arrive at
\begin{align}
\sum_{t\in\Zp}m^t\sum_{x\in\Zd}x_1^{2j}\varphi_t^\rw(x)=(2j)!\,
 \frac{(mv_\alpha)^j}{(1-m)^{j+1}}+O\big((1-m)^{-j}\big).
\end{align}
However, by the general binomial expansion,
\begin{align}\lbeq{gen-binom}
\frac{m^j}{(1-m)^{j+1}}=m^j\sum_{l=0}^\infty\binom{-j-1}l(-m)^l=m^j\sum_{l=
 0}^\infty\binom{j+l}jm^l=\sum_{t=j}^\infty\binom{t}jm^t.
\end{align}
Therefore,
\begin{align}
\sum_{x\in\Zd}x_1^{2j}\varphi_t^\rw(x)\underset{t\uparrow\infty}\sim(2j)!\,
 \binom{t}j\,v_\alpha^j\sim\frac{\Gamma(2j+1)}{\Gamma(j+1)}\,(v_\alpha t)^j.
\end{align}
This completes the proof of \refeq{main2} for $r=2j<\alpha$.

In order to consider the other values of $r<\alpha$, we use the following
integral representation for $|x_1|^q$ with $q\in(0,2)$ (cf., \cite{csIII}):
\begin{align}\lbeq{fracmom}
|x_1|^q=\frac1{K_q}\int_0^\infty\frac{1-\cos(ux_1)}{u^{1+q}}\,\text{d}u,
\end{align}
where
\begin{align}
K_q=\int_0^\infty\frac{1-\cos u}{u^{1+q}}\,\text{d}u=\frac\pi{2\sin\frac{q\pi}
 2}\,\frac1{\Gamma(q+1)}.
\end{align}
Let $r=2j+q$ with $j\in\Zp$ and $q\in(0,2)$.  Then, by \refeq{fracmom},
the generating function for the fractional moment
$\big\{\sum_{x\in\Zd}|x_1|^{2j+q}\varphi_t^\rw(x)\big\}_{t\in\Zp}$ can be
written as
\begin{align}\lbeq{rw-gen}
\sum_{t\in\Zp}m^t\sum_{x\in\Zd}|x_1|^{2j+q}\varphi_t^\rw(x)&=\frac1{K_q}
 \int_0^\infty\frac{\text{d}u}{u^{1+q}}\sum_{t\in\Zp}m^t\sum_{x\in\Zd}\big(1
 -\cos(ux_1)\big)\,x_1^{2j}\varphi_t^\rw(x)\nn\\
&=\frac{(-1)^j}{K_q}\int_0^\infty\frac{\text{d}u}{u^{1+q}}\,\Big(\nablone^{2j}
 \hat\varphi_m^\rw(0)-\nablone^{2j}\hat\varphi_m^\rw(\vec u)\Big),
\end{align}
where $\vec u=(u,0,\dots,0)\in\Rd$.  Therefore, similarly to the above case of 
$r=2j$, it suffices to investigate the ``derivative"
\begin{align}
\laplace_{\vec u}\nablone^{2j}\hat\varphi_m^\rw(0)\equiv\nablone^{2j}\hat
 \varphi_m^\rw(0)-\nablone^{2j}\hat\varphi_m^\rw(\vec u).
\end{align}
However, by ``differentiating" both sides of \refeq{rw-exp} and using the
$\Zd$-symmetry, we obtain
\begin{align}
\laplace_{\vec u}\nablone^{2j}\hat\varphi_m^\rw(0)&=m\laplace_{\vec u}
 \nablone^{2j}\hat\varphi_m^\rw(0)+m\sum_{l=1}^j\binom{2j}{2l}\nablone^{2l}
 \hat D(0)\,\laplace_{\vec u}\nablone^{2(j-l)}\hat\varphi_m^\rw(0)\nn\\
&\quad+m\sum_{n=0}^{2j}\binom{2j}n\nablone^{2j-n}\hat\varphi_m^\rw(\vec u)\,
 \laplace_{\vec u}\nablone^n\hat D(0)\nn\\
&=\frac{m}{1-m}\bigg(\sum_{l=1}^j\binom{2j}{2l}\nablone^{2l}\hat D(0)\,
 \laplace_{\vec u}\nablone^{2(j-l)}\hat\varphi_m^\rw(0)\nn\\
&\hskip4pc+\sum_{n=0}^{2j}\binom{2j}n\nablone^{2j-n}\hat\varphi_m^\rw(\vec u)\,
 \laplace_{\vec u}\nablone^n\hat D(0)\bigg),
\end{align}
where we regard the sum over $l\in\{1,\dots,j\}$ in the last expression as
zero when $j=0$.  Substituting this back to \refeq{rw-gen}, performing the
integration with respect to $u\in(0,\infty)$ and then reorganizing the
resulting terms (see \cite{csIII} for more details), we will end up with
\begin{align}
\sum_{t\in\Zp}m^t\sum_{x\in\Zd}|x_1|^r\varphi_t^\rw(x)&=\frac{2\sin\frac{r\pi}
 {\alpha\vee2}}{(\alpha\wedge2)\sin\frac{r\pi}\alpha}\,\Gamma(r+1)\,\frac{(m
 v_\alpha)^{\frac{r}{\alpha\wedge2}}}{(1-m)^{1+\frac{r}{\alpha\wedge2}}}\nn\\
&\quad\times
 \begin{cases}
 1+O((1-m)^\epsilon)&(\alpha\ne2),\\[5pt]
 \big(\log\frac1{\sqrt{1-m}}\big)^{r/2}+O(1)\quad&(\alpha=2),
 \end{cases}
\end{align}
for some $\epsilon>0$.  The proof of \refeq{main2} is completed by expanding
the right-hand side of the above expression in powers of $m$ and comparing the
coefficient of $m^t$ in both sides, for large $t$.

\section{The model-dependence}\label{s:dc}
The key to the proof for self-avoiding walk and oriented percolation is the
following lace expansion (see, e.g., \cite{csI,s06}):
\begin{align}\lbeq{lace-exp-x}
\varphi_t(x)=I_t(x)+\sum_{s=1}^t(J_s*\varphi_{t-s})(x),
\end{align}
where
\begin{align}\lbeq{IJ-def-x}
I_t(x)=
 \begin{cases}
 \delta_{x,o}\delta_{t,0}&(\text{SAW}),\\
 \pi_t^\op(x)&(\text{OP}),
 \end{cases}&&
J_t(x)=
 \begin{cases}
 D(x)\delta_{t,1}+\pi_t^\saw(x)&(\text{SAW}),\\
 p(D*\pi_{t-1}^\op)(x)&(\text{OP}).
 \end{cases}
\end{align}
Recall \refeq{rw-exp-x} for random walk, so that
$I_t^\rw(x)=\delta_{x,o}\delta_{t,0}$ and $J_t^\rw(x)=D(x)\delta_{t,1}$.  The
model-dependent $\pi_t(x)$ in \refeq{IJ-def-x}, which accounts for difference
from random walk, is an alternating sum of the lace-expansion coefficients and
obey the following diagrammatic bounds (cf., \cite{csI,s06}):
\begin{align}
&|\pi_t^\saw(x)|\le
 ~\mathop{\raisebox{-9pt}{\includegraphics[scale=.09]{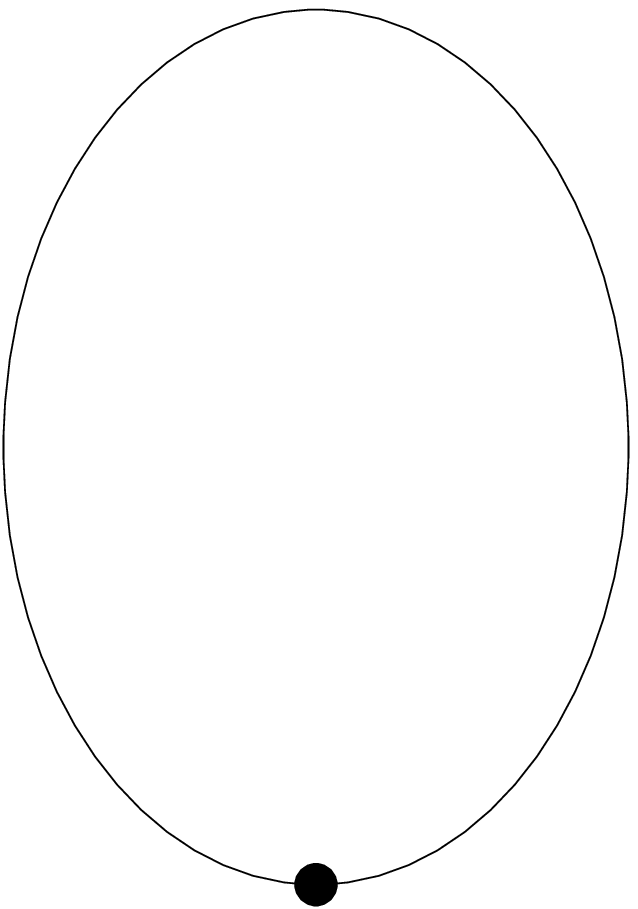}}}_{x=o}~
 +~\mathop{\raisebox{-9pt}{\includegraphics[scale=.09]{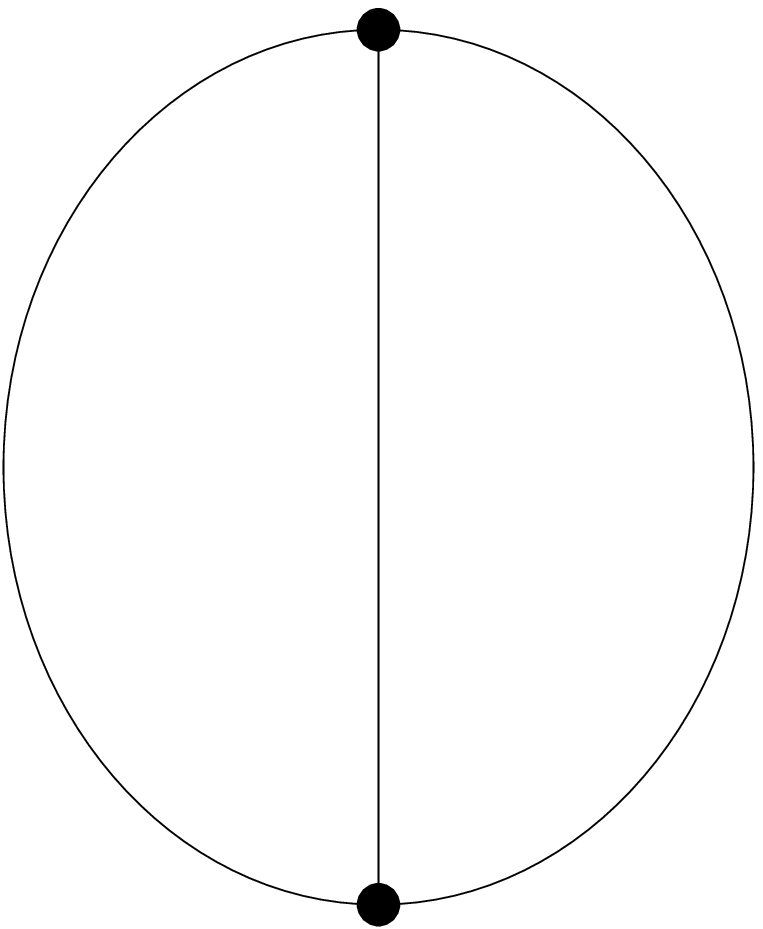}}}_o^x~
 +~\mathop{\raisebox{-9pt}{\includegraphics[scale=.09]{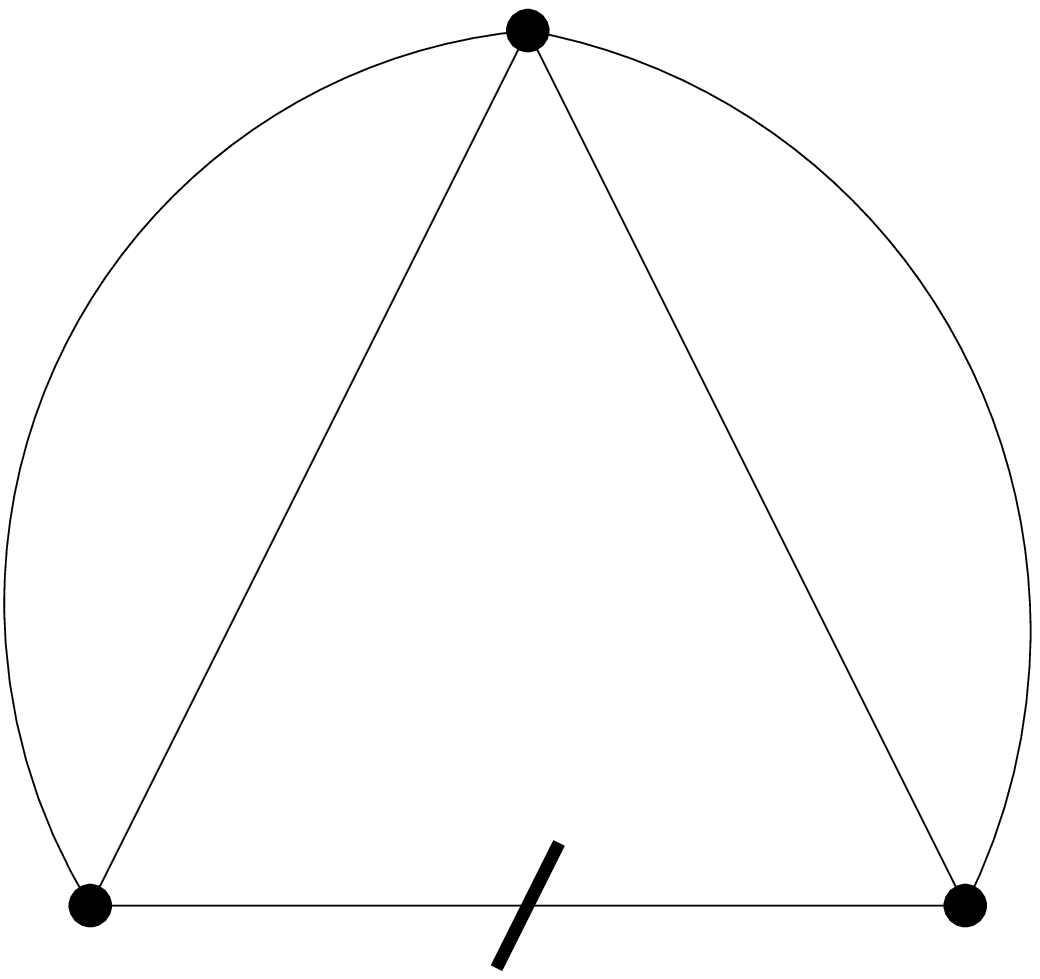}}}_{o\qquad x}~
 +\cdots,\lbeq{saw-diagbd}\\
&|\pi_t^\op(x)|\le
 ~\mathop{\raisebox{-12pt}{\includegraphics[scale=.09]{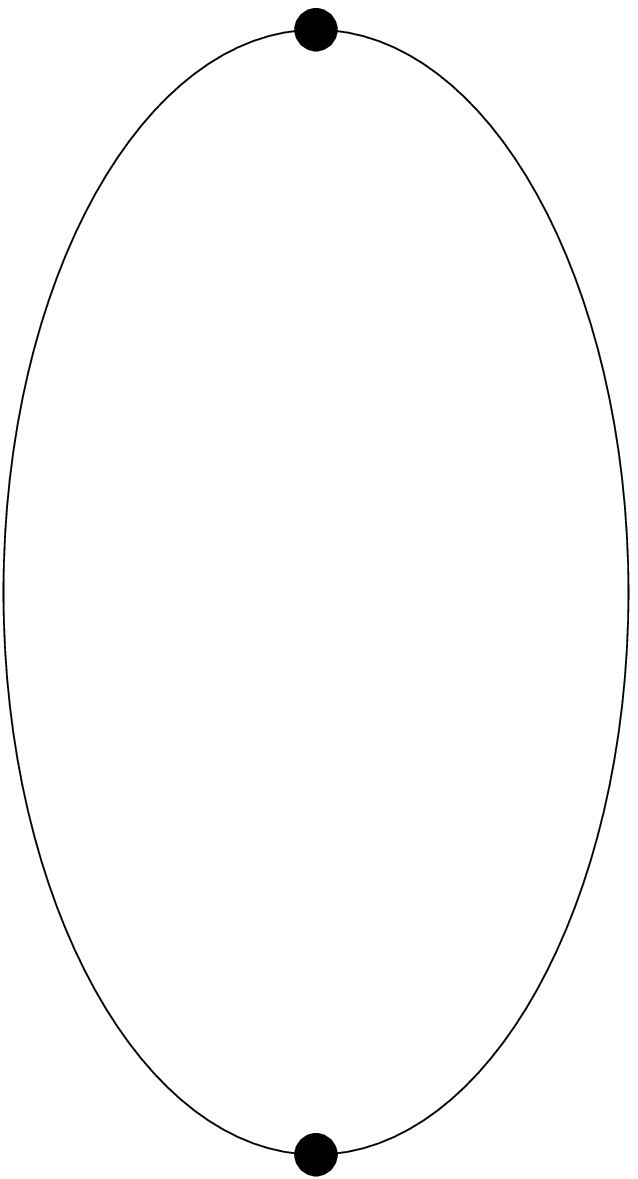}}}_{(o,0)}^{(x,t)}~
 +~\mathop{\raisebox{-14pt}{\includegraphics[scale=.1]{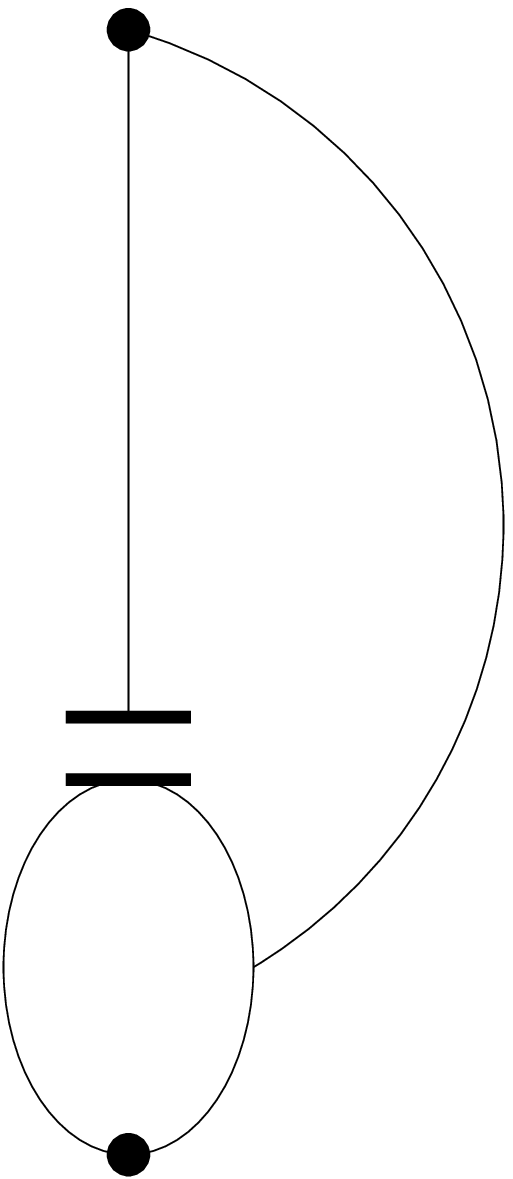}}}_{(o,0)}^{(x,t)}~
 +~\mathop{\raisebox{-17pt}{\includegraphics[scale=.12]{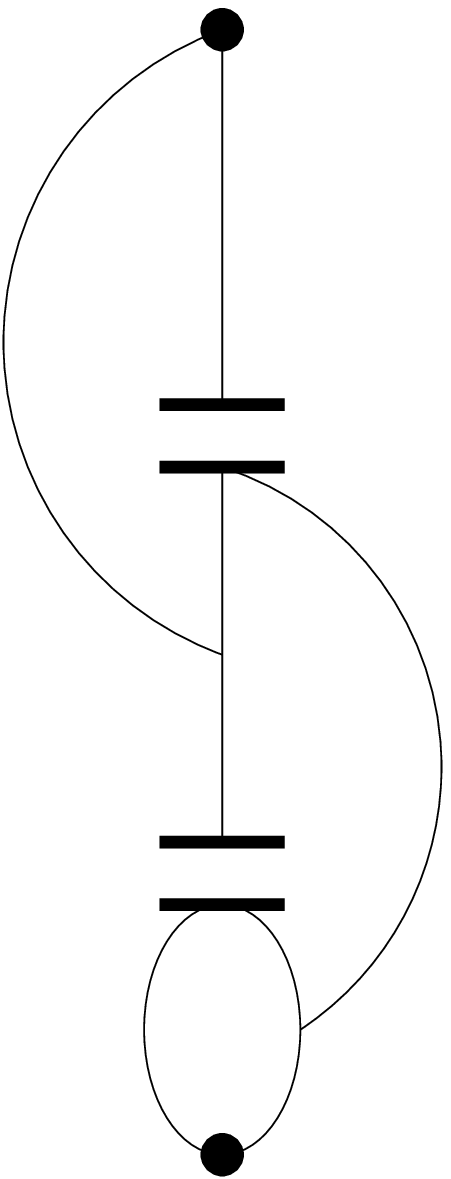}}}_{(o,0)}^{(x,t)}~
 +~\mathop{\raisebox{-17pt}{\includegraphics[scale=.12]{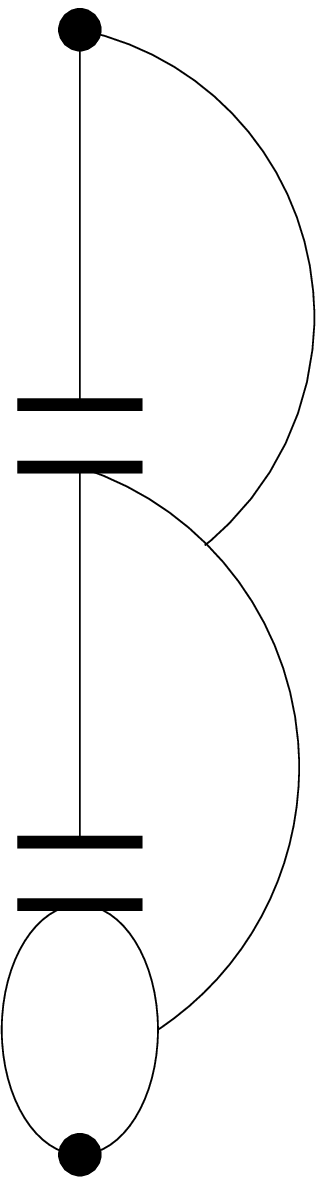}}}_{(o,0)}^{(x,t)}~
 +~\cdots,\lbeq{op-diagbd}
\end{align}
where each line corresponds to a 2-point function.  For self-avoiding walk,
the first diagram represents self-avoiding loop of length $t\ge2$, i.e.,
$(D*\varphi_{t-1}^\saw)(x)$, and the second diagram represents the product of
three 2-point functions,
$\varphi_{s_1}^\saw(x)\,\varphi_{s_2}^\saw(x)\,\varphi_{s_3}^\saw(x)$,
summed over all possible combinations of $s_1,s_2,s_3\ge1$ satisfying
$s_1+s_2+s_3=t$, and so on.  For oriented percolation, the first diagram
represents $\varphi_t^\op(x)^2$, where the upward direction is the
time-increasing direction, and the second diagram represents the product of
five 2-point functions concatenated in a depicted way, where unlabeled
vertices are summed over $\Zd\times\Zp$, and so on.  For more details,
we refer to \cite{csI,s06}.

Because of the similarity between \refeq{rw-exp-x} and \refeq{lace-exp-x},
it is natural to expect that the strategy in \S\ref{s:rw} for random walk
may also work for self-avoiding walk and oriented percolation.  To see if
it really works, we first take the Fourier-Laplace transform of
\refeq{lace-exp-x}.  For $k\in[-\pi,\pi]^d$ and $m\in[0,\mc)$,
\begin{align}
\hat\varphi_m(k)=\hat I_m(k)+\hat J_m(k)\,\hat\varphi_m(k),
\end{align}
where $\mc\ge1$ is the model-dependent radius of convergence for
$\big\{\sum_{x\in\Zd}\varphi_t(x)\big\}_{t\in\Zp}$ for self-avoiding walk
and critcal/subcritical oriented percolation ($\mc^\op$ is a non-increasing
function of $p\le\pc$ and $\mc^\op=1$ at $p=\pc$ \cite{csI}).
Due to the diagrammatic bounds \refeq{saw-diagbd}--\refeq{op-diagbd}, it has been proved
\cite{csI,csII,h??} that, for $d>\dc$ and $L\gg1$, there are
$\epsilon,\delta>0$ such that
\begin{align}
\sum_{t\in\Zd}t^{1+\epsilon}m^t\sum_{x\in\Zd}|\pi_t(x)|,&&
\sum_{t\in\Zd}m^t\sum_{x\in\Zd}|x_1|^{\alpha\wedge2+\delta}|\pi_t(x)|,
\end{align}
both converge, even at $m=\mc$.  This implies that
$\hat J_{\mc}(0)=1$ and, as $m\uparrow\mc$,
\begin{align}\lbeq{lace-exp}
\hat\varphi_m(0)=\frac{\hat I_m(0)}{1-\hat J_m(0)}=\frac{\hat I_m(0)}{\hat
 J_{\mc}(0)-\hat J_m(0)}&\sim\frac{\hat I_{\mc}(0)}{\mc\partial_m\hat J_{\mc}
 (0)\,\big(1-\frac{m}{\mc}\big)}\nn\\
&=\frac{\hat I_{\mc}(0)}{\mc\partial_m\hat J_{\mc}(0)}\sum_{t\in\Zp}\Big(
 \frac{m}{\mc}\Big)^t.
\end{align}
On the other hand, for $r=2j<\alpha$ with $j\in\N$,
\begin{align}
\nablone^{2j}\hat\varphi_m(0)&=\nablone^{2j}\hat I_m(0)+\sum_{l=0}^j\binom{2j}
 {2l}\nablone^{2l}\hat J_m(0)\,\nablone^{2(j-l)}\hat\varphi_m(0)\nn\\
&=\frac1{1-\hat J_m(0)}\bigg(\nablone^{2j}\hat I_m(0)+\sum_{l=1}^j\binom{2j}
 {2l}\nablone^{2l}\hat J_m(0)\,\nablone^{2(j-l)}\hat\varphi_m(0)\bigg).
\end{align}
Suppose that the leading contribution is due to the $l=1$ term (this is far
from trivial and needs to be proved, as in \cite{csIII}).  Then, by
induction and using \refeq{lace-exp},
\begin{align}
\nablone^{2j}\hat\varphi_m(0)&\sim\binom{2j}2\frac{\nablone^2\hat J_m(0)}{1-
 \hat J_m(0)}\,\nablone^{2(j-1)}\hat\varphi_m(0)\nn\\
&~\;\vdots\nn\\
&\sim\frac{(2j)!}{2^j}\bigg(\frac{\nablone^2\hat J_m(0)}{1-\hat J_m(0)}\bigg)^j
 \hat\varphi_m(0)\nn\\
&\sim\frac{(2j)!}{2^j}\bigg(\frac{\nablone^2\hat J_{\mc}(0)}{\mc\partial_m\hat
 J_{\mc}(0)\,\big(1-\frac{m}{\mc}\big)}\bigg)^j\frac{\hat I_{\mc}(0)}{\mc
 \partial_m\hat J_{\mc}(0)\,\big(1-\frac{m}{\mc}\big)}.
\end{align}
However, similarly to \refeq{gen-binom},
\begin{align}
\Big(1-\frac{m}{\mc}\Big)^{-j-1}=\sum_{t\in\Zp}\binom{t+j}j\Big(\frac{m}{\mc}
 \Big)^t,
\end{align}
hence
\begin{align}\lbeq{num}
\nablone^{2j}\hat\varphi_m(0)\sim(2j)!\,\bigg(\frac{\nablone^2\hat J_{\mc}(0)}
 {2\mc\partial_m\hat J_{\mc}(0)}\bigg)^j\frac{\hat I_{\mc}(0)}{\mc\partial_m
 \hat J_{\mc}(0)}\sum_{t\in\Zp}\binom{t+j}j\Big(\frac{m}{\mc}\Big)^t.
\end{align}
Therefore, by \refeq{lace-exp} and \refeq{num},
\begin{align}\lbeq{lace-exp-res}
\frac{\sum_{x\in\Zd}x_1^{2j}\varphi_t(x)}{\sum_{x\in\Zd}\varphi_t(x)}&\sim
 \frac{(2j)!}{j!}\,\bigg(\frac{-\nablone^2\hat J_{\mc}(0)}{2\mc\partial_m\hat
 J_{\mc}(0)}\,t\bigg)^j\nn\\
&=\frac{\Gamma(2j+1)}{\Gamma(j+1)}\,\bigg(\underbrace{\frac1{\mc\partial_m
 \hat J_{\mc}(0)}\,\frac{\nablone^2\hat J_{\mc}(0)}{\nablone^2\hat D
 (0)}}_{C_\alpha}\;\underbrace{\frac{-\nablone^2\hat D(0)}2}_{v_\alpha}\;t
 \bigg)^j.
\end{align}
This completes a sketch proof for $r=2j$.

The case for the other values of $r<\alpha$ is more involved, but can be
proved by following the same strategy as in \S\ref{s:rw} for random walk.
However, since $C_\alpha$ in \refeq{lace-exp-res} is ill-defined for
$\alpha\le2$ due to the divergence of $\nablone^2\hat D(0)$, it is replaced by
\begin{align}
C_\alpha=\frac1{\mc\partial_m\hat J_{\mc}(0)}\lim_{k\to0}\frac{\laplace_k\hat
 J_{\mc}(0)}{\laplace_k\hat D(0)}\equiv\frac1{\mc\partial_m\hat J_{\mc}(0)}
 \lim_{k\to0}\frac{\hat J_{\mc}(0)-\hat J_{\mc}(k)}{\hat D(0)-\hat D(k)}.
\end{align}
We refrain from showing further details and refer the readers to the original
paper \cite{csIII}.

\section*{Acknowledgements}
This work was supported by the start-up fund of the Leader Development System
in the Basic Interdisciplinary Research Areas at Hokkaido University.  I am
grateful to Lung-Chi Chen for the fruitful collaboration on the long-range
models \cite{csI,csII,csIII} and Keiichi R.~Ito for the invitation to the RIMS
workshop ``Applications of RG Methods in Mathematical Sciences" at Kyoto
University from September 9 through $11^\text{th}$, 2009.


\begin{thebibliography}{99}

\bibitem{csI}L.-C. Chen and A. Sakai.
\newblock Critical behavior and the limit distribution for long-range oriented
percolation.~I.
\newblock\emph{Probab. Theory Relat. Fields}~\textbf{142} (2008): 151--188.

\bibitem{csII}L.-C. Chen and A. Sakai.
\newblock Critical behavior and the limit distribution for long-range oriented
percolation.~II: Spatial correlation.
\newblock\emph{Probab. Theory Relat. Fields}~\textbf{145} (2009): 435--458.

\bibitem{csIII}L.-C. Chen and A. Sakai.
\newblock Asymptotic behavior of the gyration radius for long-range self-avoiding walk
and long-range oriented percolation.
\newblock In preparation.

\bibitem{h??}M. Heydenreich.
\newblock Long-range self-avoiding walk converges to $\alpha$-stable processes.
\newblock Preprint, arXiv:0809.4333v1 (2008).

\bibitem{s06}G. Slade.
\newblock The lace expansion and its applications.
\newblock\emph{Lecture Notes in Math.}~\textbf{1879} (2006).

\end{thebibliography}
\end{document}